\documentclass [12pt]  {article} 
\usepackage{amssymb}
\def \CC{{\mathbb{C}}}
\def \FFF{{\cal{F}}}
\def \HHH{{\cal{H}}}

\begin{document}

\begin{center}
{\Large {\bf A normality criterion corresponding to the defect relations}}\\
\bigskip
{\sc Andreas Schweizer}\\
\bigskip
{\small {\rm Department of Mathematics,\\
Korea Advanced Institute of Science and Technology (KAIST),\\ 
Daejeon 305-701\\
South Korea\\
e-mail: schweizer@kaist.ac.kr}}
\end{center}
\begin{abstract}
\noindent
Let $\FFF$ be a family of meromorphic functions on 
a domain $D$. We present a quite general sufficient 
condition for $\FFF$ to be a normal family. This 
criterion contains many known results as special 
cases. The overall idea is that certain comparatively 
weak conditions on $\FFF$ locally lead to somewhat 
stronger conditions, which in turn lead to even 
stronger conditions on the limit function $g$ in 
the famous Zalcman Lemma. Ultimately, the defect 
relations for $g$ force normality of $\FFF$.
\end{abstract}

\subsection*{1. Introduction and main result}

\noindent
A family $\FFF$ of meromorphic functions on a domain 
$D\subseteq\CC$ is {\bf normal} if every sequence 
$\{f_n\}\subseteq\FFF$ contains a subsequence which
converges spherically uniformly on compact subsets 
of $D$.
\par
We refer to [10] for all necessary background information.
\\ 
The starting point for this paper was the following recent result.
\\ \\
{\bf Theorem 1.1.} [16, Theorem 1.6] \it
Let $\FFF$ be a family of meromorphic functions defined in a domain 
$D$, $M$ be a positive number and $S=\{\alpha, \beta\}$, where 
$\alpha$, $\beta$ are distinct elements of $\CC\cup\{\infty\}$.
Further suppose that 
\begin{itemize}
\item[(i)] each pair of functions $f$, $g\in\FFF$ share the set
$S$ in $D$;
\item[(ii)] there exists a $\gamma\in\CC\setminus\{\alpha,\beta\}$
such that for each $f\in\FFF$, $|f'(z)|\leq M$ whenever $f(z)=\gamma$
in $D$;
\item[(iii)] each $f\in\FFF$ has no simple $\beta$-points in $D$.
\end{itemize}
Then $\FFF$ is normal in $D$.
\rm
\\ \\
Condition (iii) is a bit of a nuisance because it destroys 
the symmetry between $\alpha$ and $\beta$. And indeed, 
condition (iii) can be completely omitted. This can be seen 
as follows.
\par
To prove normality at a point $z_0$ we can restrict $\FFF$ to 
a small neighbourhood of $z_0$. Given a sequence of functions 
$f_n$ from $\FFF$, we can make this neighbourhood so small that
the sequence $f_n$ contains a subsequence $f_{n_k}$ such that in
this small neighbourhood all $f_{n_k}$ share $\alpha$ and omit
$\beta$. Or they all share $\beta$ and omit $\alpha$, in which 
case we simply switch the names.
See Step 1 in the proof of Theorem 1.2 given in Section 3 for 
more details and a slightly more general version. 
So in final instance, Theorem 1.1 itself, applied locally, 
proves the version without condition (iii). 
\par
Besides Theorem 1.1 the improved version also improves
[13, Theorem 1]. Actually, both of them, and several other
published results, can also be obtained as special cases of
our quite general main result, which we will formulate as 
Theorem 1.2.
\par
In addition to the local argument that we just outlined,
the proof of Theorem 1.2 also makes heavy use of the impact
that certain conditions satisfied by all $f\in\FFF$ will 
have on the limit function $g$ in Zalcman's Lemma. 
These steps follow the ideas from the proofs of 
[3, Theorem 2] and [8, Theorem 1.3]. Then in the
end everything will boil down to the defect relations 
applied to $g$.
\\ \\
{\bf Theorem 1.2. (Main Theorem)} \it
Let $\FFF$ be a family of meromorphic functions on 
a domain $D$. Let there be pairwise disjoint, finite 
sets $A\subset\CC\cup\{\infty\}$, $B\subset\CC^*$ and 
$C\subset\CC\cup\{\infty\}$ of cardinalities $n$, $s$ 
and $r$, and with the following properties: 
\begin{itemize}
\item[(a)] All $f\in\FFF$ share the set $A$, i.e., for every 
$z_1 \in D$ we have $f(z_1)\in A$ for all $f\in\FFF$ or 
$f(z_1)\not\in A$ for all $f\in\FFF$.
\item[(b)] For every $b_i \in B$, every $f\in\FFF$ and every 
$z_1 \in D$ we have $f(z_1)=b_i\Rightarrow f'(z_1) =b_i$.
\item[(c)] If $\infty\in C$, say $c_1 =\infty$, then on $D$ every
pole of every $f\in\FFF$ has at least multiplicity $m_1\geq 2$.
For every finite $c_j \in C$ there exists an integer 
$m_j \geq 2$ and a number $M_j$ such that for every $f\in\FFF$ 
and every $z_1 \in D$ we have 
$f(z_1)=c_j \Rightarrow |f^{(k)}(z_1)|\leq M_j$ 
for $k=1,\ldots, m_j -1$.
\end{itemize}
If
$$n+s+\sum_{j=1}^r(1-\frac{1}{m_j})>2,$$
then $\FFF$ is normal on $D$.
\rm
\\ \\
The proof will be given in Section 3.
\par
In practice the conditions will often come in a less technical 
form. If for example on $D$ all $f\in\FFF$ omit the value 
$\omega\in\CC^*$, the we could take $\omega$ into any one of 
the three sets $A$, $B$ or $C$ (with $m=\infty$), but of course
only into one because of the pairwise disjointness.
\par
A value $\omega$ such that all $\omega$-points of all
$f\in\FFF$ have multiplicity at least $m$ clearly belongs
into the set $C$ with $m_i =m$. In fact, the condition for
$C$ is weaker. 
\par
Note however that if for all $f\in\FFF$ we have 
$f\in\{d_1,\ldots,d_t\}\Rightarrow f'\in\{d_1,\ldots,d_t\}$
with $t>1$, then the values $d_1,\ldots,d_t$ are only booked
as belonging to the set $C$ (with $m=2$), even if all $d_i$ 
are nonzero. See Section 4 for more on this.
\par
Also note that from a condition $f=a\Leftrightarrow f'=a$ we 
only use one direction. For more on this, see Section 4 as well.
\\ \\
{\bf Corollary 1.3.} [14, Theorem 2] \it
Let $a_1$, $a_2$, $a_3$ be three distinct complex numbers. 
If for each $f\in\FFF$ we have $f=a_i \Rightarrow f'=a_i$
($i=1,2,3$) in $D$, then $\FFF$ is normal.
\rm
\\ \\
{\bf Proof.} \rm
Obviously at least two of the $a_i$ are nonzero. So we have 
$|B|=2$ and $|C|=1$ with $m_1 =2$.
\hfill$\Box$
\\ \\
Other corollaries, which we don't want to all spell out in detail, 
include our improved version of Theorem 1.1, [3, Theorem 2], 
[8, Theorem 1.3], [10, Theorem 4.1.4, p.105],
[11, Theorem 1.3], [13, Theorems 1 and 2], and 
[17, Theorem 5].
\par
Most of them are generalizations of Montel's Fundamental Normality
Criterion (see for example [10, p.74]), which says that a family
whose functions all omit the same three values must be normal. 
They cover many of the possibilities with $|A|+|B|+|C|=3$.
In Section 2 we discuss the cases that do not imply normality.
We present one more constellation that we couldn't find in the 
literature. 
\\ \\
{\bf Corollary 1.4.} \it
Let $\FFF$ be a family of meromorphic functions on a domain 
$D$ with the following properties:
\begin{itemize}
\item[(a)] For any $f\in\FFF$ and any $z_1 \in D$ we have 
$f'(z_1)=1$ whenever $f(z_1)=1$.
\item[(b)] All poles of all $f\in\FFF$ have multiplicity 
at least $3$.
\item[(c)] The set $\bigcup\limits_{f\in\FFF}f'(f^{-1}\{0\})$ 
is bounded.
\end{itemize}
Then $\FFF$ is normal on $D$.
\rm
\\ \\
{\bf Proof.} \rm
In Theorem 1.2 this gives $B=\{1\}$ and $C=\{\infty, 0\}$ with 
$m_1 =3$ and $m_2 =2$.
\hfill$\Box$
\\

\subsection*{2. (Counter-)Examples}

\noindent
In this section we will construct several examples that show 
that the conditions in Theorem 1.2 are sharp. What we mean by 
that is that for numbers $n$, $s$, $r$ and $m_1 ,\ldots, m_r$ 
such that $n+s+\sum_{k=1}^r(1-\frac{1}{m_k})\leq 2$ one can 
construct a non-normal family satisfying the other conditions 
of Theorem 1.2 for these numbers. 
\par
We need some preparation. The spherical derivative of a meromorphic 
function $f$ is 
$$f^{\#}(z):=\frac{|f'(z)|}{1+|f(z)|^2}.$$

\bigskip
\noindent
{\bf Theorem 2.1.} [2, Lemma 1] \it
Let $f$ be a meromorphic function on $\CC$. If $f$ has bounded 
spherical derivative on $\CC$, $f$ is of order at most $2$. If,
in addition, $f$ is entire, then the order of $f$ is at most $1$.
\rm
\\ \\
{\bf Corollary 2.2.} \it
Let $f$ be a meromorphic function on $\CC$. If the order of $f$ 
is bigger than $2$, the family
$$\FFF=\{f(z+\omega)\ :\ \omega\in\CC\}$$ 
is not normal at any point $z_0 \in\CC$.
\rm
\\ \\
{\bf Proof.} \rm 
By Marty's Criterion [10, p.75], a family $\FFF$ of meromorphic 
functions on a domain $D$ is normal if and only if for every compact
subset $K\subseteq D$ there exists a constant $C(K)$ such that
$f^{\#}(z)\leq C(K)$ for all $z\in K$ and all $f\in\FFF$. With this, 
the corollary is an immediate consequence of Theorem 2.1.
\hfill$\Box$
\\ \\
One of the many uses of normal families is the following.
Given a function $f$ that is meromorphic in the complex plane,
and certain sharing conditions for $f$ and $f'$, one wants to
show $f=f'$ or that $f$ has a certain form. In some cases 
a good strategy is not to start immediately with arguments 
from Nevanlinna Theory, but to first look at the family 
$\{f(z+\omega)\ :\ \omega\in\CC\}$.
If this family turns out to be normal, this implies that the 
order of the function $f$ is at most $2$, which might be quite 
helpful for the successive arguments (for example from Nevanlinna 
Theory). This is for example the approach in [8].
\par
Here we will go exactly in the opposite direction. We find 
a function $h$ of order bigger than $2$ that has some desired 
properties like omitting a value, or sharing a value with its 
derivative, or having a totally ramified value. Because of the 
order of $h$ the family $\{h(z+\omega)\ :\ \omega\in\CC\}$ cannot
be normal. But it inherits the desired conditions. This proves 
that these conditions are not sufficient for normality.
\par
Constructing suitable functions of order bigger than $2$ is 
actually quite easy, using the following result.
\\ \\
{\bf Theorem 2.3.} [5, Corollary 1.2] \it
Let $F$ be a meromorphic function which is not of order zero 
and let $g$ be an entire function which is not a polynomial.
Then $F(g(z))$ is of infinite order.
\rm
\\ \\
{\bf Example 2.4.} 
Consider the function 
$$h(z)=\sin(e^z).$$ 
It has infinite order by Theorem 2.3 as $\sin(z)$ 
has order one. So by Corollary 2.2 the family 
$$\FFF=\{h(z+\omega)\ :\ \omega\in\CC\}$$ 
cannot be normal. Obviously, every $f\in\FFF$ omits the value 
$\infty$. Moreover, every $1$-point and every $-1$-point of 
every $f\in\FFF$ has multiplicity $2$. So in the sense of 
Theorem 1.2 we have $A=\{\infty\}$, $B=\emptyset$ and 
$C=\{1,-1\}$ with $m_1 =m_2 =2$. 
\par
Applying a linear transformation to get all three interesting
values into $\CC^*$, e.g. considering
$$h(z)=2-\frac{1}{\sin(e^z)},$$
this can also be interpreted as an example for $A=\emptyset$,
$B=\{2\}$, $C=\{1,3\}$ with $m_1 =m_2 =2$, or alternatively as 
an example for $A=B=\emptyset$, $|C|=3$ with $m_1 =m_3 =2$ and
$m_2\geq 2$.
\par
We can also apply a linear transformation to place the totally
ramified values at $\infty$ and $0$. For example, 
$$h(z)=\frac{\sin(e^z)-1}{\sin(e^z)+1}$$
is of infinite order, has only multiple zeroes and multiple 
poles, and omits the value $1$. So the non-normal family 
$\{h(z+\omega)\ :\ \omega\in\CC\}$ 
answers the question asked in [12, Remark 5]. Even 
if the function $\psi$ in [12, Theorem 2] is constant $1$, the 
condition $(2)$ in that theorem is best possible in the sense 
that it cannot be weakened to ``all poles of $f$ are multiple''.
\\ \\
{\bf Example 2.5.} 
This example is included to show that in Theorem 1.2 the 
condition $f=0\Rightarrow f'=0$ is really weaker than 
$f=b_i \Rightarrow f'=b_i$ for a nonzero complex number 
$b_i$. Let $\FFF$ be a family of analytic functions. 
If for all $f\in\FFF$ we have $f=1\Rightarrow f'=1$, and 
all $c$-points of $f$ are multiple, then Theorem 1.2 
implies that $\FFF$ is normal. On the other hand, the 
analytic family 
$\FFF=\{h(z+\omega)\ :\ \omega\in\CC\}$ 
with 
$$h(z)=\frac{c}{2}(\sin(e^z)+1),$$ 
where $c$ is a fixed nonzero complex number, is not normal. 
But for every $f\in\FFF$ we have $f=0\Rightarrow f'=0$, 
and every $c$-point of $f$ is multiple.
\\ \\
{\bf Example 2.6.} 
The function $h(z)=\exp(z^3)$ has order $3$. So by 
Corollary 2.2 the family 
$\FFF=\{h(z+\omega)\ :\ \omega\in\CC\}$
is not normal. In the terminology of Theorem 1.2 this 
is an example with $A=\{\infty,0\}$, $B=C=\emptyset$
or with $A=\{\infty\}$, $B=\emptyset$, $C=\{0\}$ and
$m_1$ any integer we want. 
\par
Moving the two omitted values to $a$ and $b$ by considering
$$h(z)=\frac{a e^{z^3}+b}{e^{z^3}+1}$$
also furnishes examples with $A=B=\emptyset$, $C=\{a,b\}$,
and $m_1$ and $m_2$ as desired, and so on. 
\par
Alternatively, we could do everything in this example
using $\exp(\exp(z))$ instead of $\exp(z^3)$. 
\\ \\
{\bf Example 2.7.} 
Fix two distinct nonzero complex numbers $a,b$ and let 
$\wp(z)$ be the Weierstrass function with differential 
equation
$$(\wp'(z))^2 =4\wp(z)(\wp(z)-a)(\wp(z)-b).$$
Since $\wp$ has order one, 
$$h(z)=\wp(e^z)$$
has infinite order by Theorem 2.3. So by Corollary 2.2 the
meromorphic family $\FFF=\{h(z+\omega)\ :\ \omega\in\CC\}$
cannot be normal. 
If $f\in\FFF$, then all poles, zeroes, $a$-points and 
$b$-points of $f$ have multiplicity $2$. This is an example 
with $A=B=\emptyset$, $|C|=4$ and $m_1=\cdots =m_4 =2$.
A similar but more complicated example is given in [7].
\par
If $b=-a$, we can also take 
$$h(z)=(\wp(e^z))^2.$$
Then for every $f\in\FFF=\{h(z+\omega)\ :\ \omega\in\CC\}$
all poles and zeroes have multiplicity $4$, and all 
$a^2$-points have multiplicity $2$. This is an example with 
$A=B=\emptyset$, $|C|=3$ and $m_1 =m_2 =4$, $m_3 =2$.
\\ \\
{\bf Example 2.8.} 
Consider the Weierstrass equation
$$(\wp'(z))^2 =4(\wp(z))^3 +c$$
with $c\neq 0$. Then $\wp'(z)$ is an elliptic function.
All its poles are triple. Differentiating the differential
equation we see that its $\sqrt{c}$-points and its 
$-\sqrt{c}$-points also all have multiplicity $3$. 
Taking 
$h(z)=\wp'(e^z)$
we get a non-normal family with 
$A=B=\emptyset$, $|C|=3$ and $m_1 =m_2 =m_3 =3$.
\par
Moreover, $(\wp'(z))^2$ has poles of order $6$, zeroes
of order $2$ and $c$-points of multiplicity $3$. So
$$h(z)=(\wp'(e^z))^2$$
gives a non-normal family with 
$A=B=\emptyset$, $|C|=3$, $m_1 =6$, $m_2 =2$, $m_3 =3$.
\\

\subsection*{3. Proof of the Main Theorem}

\noindent
Like practically every paper on normal families from the last
decade, our paper too makes heavy use of the very powerful 
Zalcman Lemma.
\\ \\
{\bf Theorem 3.1. (Zalcman's Lemma)} [15] \it
Let $\FFF$ be a family of meromorphic (analytic) functions on 
the unit disk. If $\FFF$ is not normal at $0$, then there exist
\begin{itemize}
\item[(i)] a number $0<r<1$,
\item[(ii)] points $z_n$, $|z_n|<r$,
\item[(iii)] functions $f_n\in\FFF$,
\item[(iv)] positive numbers $\rho_n\to 0$,
\end{itemize}
such that
$$f_n(z_n +\rho_n\xi)=:g_n(\xi)\to g(\xi)$$
spherically uniformly (uniformly) on compact subsets of $\CC$, 
where $g$ is a nonconstant meromorphic (entire) function.
\rm
\\ \\
The second tool that we crucially need is the following.
\\ \\
{\bf Theorem 3.2. (defect relations, slightly modified)} \it
Let $f$ be a non-constant meromorphic function on the complex plane. 
For each totally ramified value $c_j$ of $f$ let $m_j\geq 2$ be an 
integer such that each $c_j$-point of $f$ has multiplicity at least 
$m_j$. 
\begin{itemize}
\item[(a)] If $f$ is a transcendental meromorphic function, 
let $U\subset\CC\cup\{\infty\}$ be the set of all values that 
$f$ takes at only finitely many points. 
\item[(b)] If $f$ is a non-constant rational function, let 
$U\subset\CC\cup\{\infty\}$ be the set of all values that $f$ 
(as a function on the complex plane) omits. 
\item[(c)] Alternatively, if $f$ is a non-constant rational 
function and there is a value $a_1\in\CC\cup\{\infty\}$ such that 
there is excactly one $z_1\in\CC$ with $f(z_1)=a_1$, we can 
also take for $U$ the union of $\{a_1\}$ and all values that 
$f$ (as a function on the complex plane) omits. 
\end{itemize}
In any case, let $R=\{c_1,\ldots,c_r\}$ be the set of 
totally ramified values of $f$ that are not in $U$. 
Then 
$$|U|+\sum_{j=1}^r(1-\frac{1}{m_j})\leq 2.$$
\rm

\noindent
{\bf Proof.} \rm
With $U$ as in part (a) or (b) this is a special case of the 
defect relations [4, \S 5.2]. 
\par
For the proof when $U$ is as in (c), we first recall that if
$\phi:Y\to X$ is a covering of degree $d$ of compact Riemann
surfaces, then by the Hurwitz formula we have
$$2g(Y)-2=d(2g(X)-2)+\sum_{\omega\in X}(d-\#\{\phi^{-1}(\omega)\})$$
where $g$ denotes the genus. Since a rational function $f$ of 
degree $d\geq 1$ is a covering from the Riemann sphere to the
Riemann sphere, its total ramification is 
$$\sum_{\omega\in\CC\cup\{\infty\}}(d-\#\{f^{-1}(\omega)\})=2d-2.$$
Also, a rational function $f$, considered as a function on the
complex plane, can omit at most one value,
namely $f(\infty)$. In that case $f(\infty)$ must be taken 
with multiplicity $d$. Moreover, $a_1$ must then also be taken 
with multiplicity $d$. So together they use up all ramification.
Hence we have $|U|=2$ and $R$ is empty.
\par
If $f$ omits no values, then $f^{-1}(a_1)$ can consist of one 
or two points, the second one being $\infty$. In either case
the value $a_1$ uses up at least $d-2$ from the total ramification.
What is left suffices for at most two totally ramified values
(with $m_1 =m_2 =2$) or at most one totally ramified value
(possibly with bigger $m$).
\hfill$\Box$
\\ \\
The function $f(z)=\frac{(z-1)^2}{z^2 +1}$ takes each of the 
values $0$, $1$ and $2$ exactly once (on $\CC$). This shows 
that for rational functions we cannot take $U$ to be the set 
of all values that are taken at most once.
\\ \\
{\bf Proof of Theorem 1.2.} 
As normality is a local property, we can prove it at each 
point individually. So fix $z_0 \in D$.
\par
{\bf Step 1:} We pick a function $f_0$ from $\FFF$.
\par
Let us first examine the case $f_0(z_0)\in A$. If $f_0$ is 
constant, then, because of the sharing of $A$, all $f\in\FFF$ 
are constants from $A$, and thus $\FFF$ is normal. If $f_0$
is not constant, there is a small disk around $z_0$ such that
on the punctured disk $f_0$ omits all values from $A$. By 
condition (i) the same must hold for all $f\in\FFF$. We can
replace $D$ by this small disk and $\FFF$ by its restrictions
to this small disk. 
\par
Now let $f_n$ be a sequence of functions from (our new) 
$\FFF$. Our task is to show that it contains a subsequence 
that converges spherically uniformly. 
Since $f_n(z_0)\in A$ and $A$ is finite, there exists 
(at least) one element of $A$ (after renaming we can assume 
that it is $a_1$) with $f_n(z_0)=a_1$ for infinitely many $n$. 
We replace the sequence $f_n$ with the subsequence for these 
$n$. Our task still is to show that this subsequence contains 
a subsubsequence that converges spherically uniformly. For 
that we can now assume that all $f\in\FFF$ share the value 
$a_1$ and omit the values $a_2,\ldots,a_n$. The conditions 
related to the sets $B$ and $C$ remain intact.
\par
If $f_0(z_0)\not\in A$, we can make the disk around $z_0$ so 
small that $f_0$, and hence every $f\in\FFF$ completely omits 
all values from $A$. So in that case we can assume without 
loss of generality that all $f\in\FFF$ omit the set $A$
and the conditions (ii) and (iii) still hold.
\par
{\bf Step 2:} Now we assume that the family $\FFF$ with the 
stronger conditions from Step 1 is not normal at $z_0$. In
a few more steps we bring this to a contradiction. 
\par
To that end we invoke Zalcman's Lemma. Without loss of 
generality we can assume that $z_0=0$ and that $D$ is the 
unit disk.
\par
We begin by showing that the limit function $g(\xi)$ from 
Zalcman's Lemma also omits the values $a_2,\ldots,a_n \in A$
and takes the value $a_1$ at most once. 
\par
For the omitted values this follows easily from 
Hurwitz's Theorem [10, Corollary 3.8.2, p.98] 
because $g$ is not constant.
\par
If $f(0)=a_1$ for all $f\in\FFF$, then obviously $g_n(\xi)$
takes the value $a_i$ only at $\xi=\frac{-z_n}{\rho_n}$.
So, since $g$ is not constant, by Hurwitz's Theorem it can
take the value $a_1$ at a point $\xi_1$ only if every 
neighbourhood of $\xi_1$ contains $\frac{-z_n}{\rho_n}$ for
all big enough $n$. Thus 
$\lim\limits_{n\to\infty}\frac{-z_n}{\rho_n}$,
if it exists, is the only possible candidate for $\xi_1$.
\par
{\bf Step 3:} Next we show that $g$ takes each value $c_j \in C$ 
with multiplicity at least $m_j$. For fixed $n$ we differentiate 
$g_n(\xi)$ with respect to $\xi$ and get
$$g_n^{(k)}(\xi)=\rho_n^kf_n^{(k)}(z_n+\rho_n\xi).$$
Now suppose that $g(\xi_0)=c_j\in\CC$. Then there is a small 
neighbourhood of $\xi_0$ on which $g$ is holomorphic. Because 
of the locally uniform convergence, for big enough $n$ the 
functions $g_n$ must be also be holomorphic in this 
neighbourhood. Since $g$ is nonconstant, by Hurwitz's theorem 
there exist $\xi_n$, $\xi_n\to\xi_0$, such that for sufficiently 
large $n$ we have $c_j =g(\xi_0)=g_n(\xi_n)$. Hence our 
assumptions imply $|g_n^{(k)}(\xi_n)|\leq\rho_n^k M_j$ for 
$k=1,2,\ldots,m_j -1$ and $n$ sufficiently large. Since 
$g_n^{(k)}(\xi)$ converges locally uniformly to 
$g^{(k)}(\xi)$, we obtain
$$g^{(k)}(\xi_0)=\lim_{n\to\infty}g_n^{(k)}(\xi_n)=0$$
for $k=1,2,\ldots,m_j -1$. 
\par
If $c_j =\infty$, we use that 
$\frac{1}{g_n(\xi)}\to\frac{1}{g(\xi)}$ locally uniformly
(compare [10, Theorem 3.1.3, p.72]) and that all zeroes of
$\frac{1}{g_n(\xi)}$ have multiplicity at least $m_j$.
\par
{\bf Step 4:} Now we show that $g$ also omits all values from 
the set $B$. Let $b\in B$. To start with, by exactly the same 
argument as in Step 3 we see that $g$ has no simple $b$-points.
\par
So suppose $g(\xi_0)=b\in B$. Let $m$($\geq 2$) be the 
multiplicity of this $b$-point. Then $g^{(m)}(\xi_0)\neq 0$.
We can find a small disk around $\xi_0$ such that none of
$g(\xi)$, $g'(\xi), \ldots , g^{(m)}(\xi)$ vanishes at any 
point of the punctured disk. Now we fix an even smaller disk 
$\widetilde{D}$ around $\xi_0$. Since $g$ is not constant, 
by Rouch\'e's Theorem, for big enough $n$ the function $g_n$ 
must have $m$ $b$-points $\xi_{n,1},\ldots, \xi_{n,m}$ (counted
with mulitplicities) in $\widetilde{D}$. Actually, 
$$g_n'(\xi_{n,j})=\rho_n f_n'(z_n +\rho_n \xi_{n,j})=\rho_n b\neq 0;$$
so these $b$-points are all simple. Moreover, 
$\lim_{n\to\infty}g_n'(\xi_{n,j})=\lim_{n\to\infty}\rho_n b=0$.
Since $g_n'(\xi)-\rho_n b$ has $m$ zeroes on $\widetilde{D}$,
again by Rouch\'e's Theorem $g'$ must have $m$ zeroes on 
$\widetilde{D}$ (counted with multiplicity). But the only
zero of $g'$ on $\widetilde{D}$ is $\xi_0$ with multiplicity
$m-1$. This contradiction disproves the assumption that $g$
has a $b$-point.
\par
{\bf Step 5:} In Steps 2 to 4 we have shown that in the sense
of Theorem 3.2 for the non-constant function $g$ we have 
$U=A\cup B$ and $R=C$ with the same $m_j$. This yields the 
final contradiction (between the condition in Theorem 1.2 and
the conclusion of Theorem 3.2). 
\hfill$\Box$
\\

\subsection*{4. Limitations of the criterion}

\noindent
In [9] Pang and Zalcman proved, among 
other results, the following strong statement.
\\ \\
{\bf Theorem 4.1.} [9, Theorem 2] \it 
Let $\FFF$ be a family of meromorphic functions on the unit disk $D$,
and let $a$ and $b$ be distinct complex numbers. If $f$ and $f'$ share
$a$ and $b$ for every $f\in\FFF$, then $\FFF$ is normal on $D$.
\rm
\\ \\
This is a result that our approach cannot achieve. Let's even assume 
that $a$ and $b$ are both nonzero. Then the impact of the sharing 
$f=a\Leftrightarrow f'=a$ on the limit function $g$ in our proof
is that $g$ omits $a$, the same as if we only had $f=a\Rightarrow f'=a$.
And similarly for $b$. But the conditions $f=a\Rightarrow f'=a$ and
$f=b\Rightarrow f'=b$ only imply normality if $\FFF$ is an analytic
family. That's the content of [8, Theorem 1.3]. For 
a counter-example with meromorphic $f$ see Example 2.6.
\par
We also point out that if every $f\in\FFF$ shares the two-element-set
$S=\{a,b\}$ with its derivative $f'$, this does not suffice to imply 
normality, not even when the family is analytic. The following 
counter-example shows up in several papers, e.g. [3, Example 1]
and [1, Example 6]. 
\\ \\
{\bf Example 4.2.}
Consider the analytic family
$\FFF=\{f_n(z)\ :\ n=2,3,4,\ldots\}$ on the unit disk where 
$$f_n(z)=\frac{n+1}{2n}e^{nz}+\frac{n-1}{2n}e^{-nz}.$$
One easily checks that 
$n^2(f_n^2-1)=(f_n')^2 -1$.
So $f_n$ and $f_n'$ share the set $S=\{1,-1\}$. 
But $\FFF$ is not normal at $0$.
\\ \\
See also [1, Examples 7 and 8] for a non-normal meromorphic 
family in which $f$ and $f'$ share the set $\{0,b\}$ resp. the 
set $\{-1,3\}$. Actually, [1] contains a much deeper study 
of families whose functions share a two-element-set with their 
derivative.
\\ \\
{\bf Example 4.3.}
Now we consider the family $\HHH$ of all $h(z)=\frac{1}{f(z)}$
with $f(z)$ from the family $\FFF$ in Example 4.2 Then $\HHH$ 
is also not normal on the unit disk.
Since $h'(z)=\frac{f'(z)}{(f(z))^2}$, we see that the family 
$\HHH$ still satisfies the partial sharing condition 
$$h(z)\in\{1,-1\}\Rightarrow h'(z)\in\{1,-1\}.$$
Also, all zeroes are multiple, for the trivial reason that 
there are no zeroes. 
\par
This example shows that Theorem 3 in [3] is incorrect.
\\ \\
The proofs of Theorem 1 and Theorem 3 in [3] are similar, and 
the mistake seems to be that on lines 6 and 7 of page 1476 the value 
$a_l$ in $g_n'(\xi_n^{(j)})=\rho_n a_l$ depends on $j$, and therefore 
on line 13 one cannot conclude that $g_n'(\xi)-\rho_n a_l$ has $k$ 
zeroes. 
\par
So the status of Theorem 1 in [3], for which we neither know an 
alternative proof nor a counter-example, is unclear at the moment.
That theorem claims that if $S=\{a_1 ,a_2 ,a_3\}$ with three nonzero
complex numbers $a_i$, and if $f\in S\Rightarrow f'\in S$ for all
$f\in\FFF$, then $\FFF$ is normal. So apart from the condition
$a_i \neq 0$ it would be a common generalization of 
Corollary 1.3 and the following result.
\\ \\
{\bf Theorem 4.4.} [6, Theorem 1] \it 
Let $\FFF$ be a family of meromorphic functions on the unit 
disk $D$, and let $a_1$, $a_2$ and $a_3$ be three distinct 
complex numbers. If for every $f\in\FFF$, $f$ and $f'$ share 
the set $S=\{a_1 ,a_2 ,a_3 \}$, then $\FFF$ is normal on $D$.
\rm
\\ \\
The reason why Theorem 1.2 cannot yield this result is that
from the information $f\in S\Leftrightarrow f'\in S$ it only 
sees $f\in S\Rightarrow f' \in S$, and for that to force 
normality it would need that $S$ is a bounded set with at 
least $5$ elements.
\\ \\
The problematic point in the proof of Theorems 1 and 3 in
[3] is exactly the attempt to show that for a finite set 
$S$ of nonzero $a_i$ the condition $f\in S\Rightarrow f'\in S$
also forces the limit function $g$ to omit $a_i$.
\\

\subsection*{\hspace*{10.5em} References}
\begin{itemize}

\item[{[1]}] J.~Chang and Y.~Wang: \rm Shared values, Picard
values and normality, \it Tohoku Math. J. \bf 63 \rm (2011), 149-162 

\item[{[2]}] J.~Chang and L.~Zalcman: \rm Meromorphic functions
that share a set with their derivatives, \it J. Math. Anal. Appl. 
\bf 338 \rm (2008), 1020-1028 

\item[{[3]}] J.-F.~Chen: \rm Shared sets and normal families 
of meromorphic functions, \it Rocky Mountain J. Math. 
\bf vol. 40, no. 5 \rm (2010), 1473-1479

\item[{[4]}] W.~Cherry and Z.~Ye: \it Nevanlinna's Theory of
Value Distribution, \rm Springer, Berlin Heidelberg New York, 2001

\item[{[5]}] A.~Edrei and W.~H.~J.~Fuchs: \rm On the zeros of
$f(g(z))$ where $f$ and $g$ are entire functions, \it J. Analyse Math.
\bf 12 \rm (1964), 243-255 

\item[{[6]}] X.J.~Liu and X.C.~Pang: \rm Shared Values 
and Normal Families (Chinese with English Abstract), 
\it Acta. Math. Sinica \bf 50 \rm (2007), 409-412 

\item[{[7]}] F.~L\"u: \rm Corrigendum to `` A note on 
meromorphic functions that share a set with their derivatives'', 
\it Arch. Math. (Basel) \bf 97 \rm (2011), 259-260 

\item[{[8]}] F.~L\"u, J.~Xu and H.~Yi: \rm Uniqueness 
theorems and normal families of entire functions and their derivatives, 
\it Ann. Polon. Math. \bf 95.1 \rm (2009), 67-75

\item[{[9]}] X.~Pang and L.~Zalcman: \rm Normality and shared values, 
\it Ark. Mat. \bf 38 \rm (2000), 171-182

\item[{[10]}] J.~L.~Schiff: \it Normal families, 
\rm Springer, New York, 1993

\item[{[11]}] D.~Sun and H.~Liu: \rm Normal criteria 
on the family of meromorphic functions with shared set, 
\it Israel J. Math. \bf 184 \rm (2011), 403-412 

\item[{[12]}] Y.~Xu: \rm On Montel's theorem and Yang's problem,
\it J. Math. Anal. Appl. \bf 305 \rm (2005), 743-751 

\item[{[13]}] Y.~Xu: \rm Montel's criterion and shared function,
\it Publ. Math. Debrecen \bf 77/3-4 \rm (2010), 471-478 

\item[{[14]}] Y.~Xu and M.~Fang: \rm A note on some results 
of Schwick,
\it Bull. Malays. Math. Sci. Soc. \bf 27 \rm (2004), 1-8 

\item[{[15]}] L.~Zalcman: \rm A Heuristic Principle in Complex 
Function Theory, \it Amer. Math. Monthly \bf 82 \rm (1975), 813-817

\item[{[16]}] S.~Zeng and I.~Lahiri: \rm Montel's Criterion and Shared 
Set, \it Bull. Malays. Math. Sci. Soc. \bf 38 \rm (2015), 1047-1052

\item[{[17]}] Q.~Zhang: \rm Normal criteria concerning sharing values,
\it Kodai Math. J. \bf 25 \rm (2002), 8-14 

\end{itemize}

\end{document}